\documentclass[a4paper,reqno]{amsart}

\usepackage[utf8]{inputenc}
\usepackage{amssymb}
\usepackage{enumitem}
\usepackage{mathrsfs}
\usepackage{mathtools}
\usepackage{amsmath}
\usepackage{xcolor}
\usepackage[colorlinks=true]{hyperref}
\hypersetup{linkcolor=violet,urlcolor=cyan,citecolor=cyan}

\setlist[enumerate]{label=\emph{(\roman*)}}

\usepackage{environ}
\usepackage{cancel}

\usepackage{amsmath,amssymb,mathrsfs}

\newtheorem{theorem}{Theorem}[section]
\newtheorem{corollary}[theorem]{Corollary}
\newtheorem{lemma}[theorem]{Lemma}
\newtheorem{proposition}[theorem]{Proposition}

\theoremstyle{definition}

\newtheorem{remark}[theorem]{Remark}
\numberwithin{equation}{section}

\begin{document}

\parindent=0pt

	\title[1D wave-Klein-Gordon]
	{Global existence for a 1D wave--Klein--Gordon system with null coupling}
	
\author[J. Zhao]{Jingya Zhao}
\address{State Key Laboratory of Mathematical Sciences, Academy of Mathematics and Systems Science, Chinese Academy of Sciences, Beijing 100190, China}
\email{zhaojingya@amss.ac.cn}

\begin{abstract}
We prove small-data global existence for a one-dimensional wave--Klein--Gordon system in which the null form \(w_t v_x-w_x v_t\) drives both equations. For smooth, compactly supported data of
Sobolev norm at most \(\varepsilon\), the lower-order energies remain uniformly \(O(\varepsilon)\). We also establish the sharp pointwise bound
\(|v|\lesssim\varepsilon t^{-1/2}\) for the Klein--Gordon field. The proof uses a potential for the wave field to eliminate the
quadratic Klein--Gordon coupling before applying hyperboloidal energy estimates.
\end{abstract}

	\maketitle
	
\tableofcontents	
	
\section{Introduction}\label{sec:intro}
\subsection{Model problem and main result}

In this paper, we consider a wave field \(w\) and a Klein--Gordon field \(v\) of unit
mass on \(\mathbb R^{1+1}\). The two fields are coupled through the
antisymmetric null form \(Q_{01}(w,v)=w_t v_x-w_x v_t\) in the system
\begin{equation}\label{eq:wave-KG}
\begin{aligned}
    (\partial_t^2-\partial_x^2) w =Q_{01}(w,v),\\
    (\partial_t^2-\partial_x^2+1) v=Q_{01}(w,v).
\end{aligned}
\end{equation}
Here and below, coordinate subscripts denote partial derivatives, as in
\(v_t=\partial_t v\) and \(v_x=\partial_x v\). We write \(A\lesssim B\) if \(A\le CB\) for some constant \(C>0\). By time translation invariance, we prescribe the Cauchy data at \(t=2\)
and write
\begin{equation}\label{eq:intialdata}
    (w,\partial_t w)|_{t=2}=(w_0,w_1),\qquad (v,\partial_t v)|_{t=2}=(v_0,v_1).
\end{equation}


\begin{theorem}
\label{thm:main-result}
Let \(N\ge 6\) be an integer. There exists an
\(\varepsilon_0=\varepsilon_0(N)>0\) such that, if the smooth initial
data are supported in \([-1,1]\) and satisfy
\begin{equation}
  \label{eq:initial-data-smallness}
  \|v_0\|_{H^{N+2}(\mathbb R_x)}+\|w_0\|_{H^{N+2}(\mathbb R_x)}
  +\|v_1\|_{H^{N+1}(\mathbb R_x)}+\|w_1\|_{H^{N+1}(\mathbb R_x)}
  \le \varepsilon \le \varepsilon_0,
\end{equation}
then the Cauchy problem \eqref{eq:wave-KG}--\eqref{eq:intialdata} has a unique global smooth solution \((w,v)\) on
\([2,\infty)\times\mathbb R\).
The solution satisfies
\[
  |w(t,x)|+t^{1/2}|v(t,x)|\lesssim \varepsilon.
\]
\end{theorem}

\begin{remark}[Sobolev bounds on hyperbolas]
For each \(\tau\ge2\), let
\(t=\tau\cosh y\), \(x=\tau\sinh y\).
In these coordinates, the solution in Theorem~\ref{thm:main-result} satisfies
\[
  \|\partial_y w(\tau,\cdot)\|_{H_y^{N-1}}\lesssim\varepsilon,\qquad
  \|v(\tau,\cdot)\|_{H_y^{N-1}}\lesssim\varepsilon \tau^{-1/2},
\]
where \(H_y^{N-1}=H^{N-1}(\mathbb R_y)\) denotes the Sobolev space
in the variable \(y\) at fixed \(\tau\).

The first bound follows from \eqref{eq:improved-low}, since
\(\partial_y w\) is the second component of \(U\) defined in
Section~\ref{sec:potential-system}. The second follows from
\eqref{eq:improved-low} and Lemma~\ref{lem:v-sobolev} with \(m=N-1\).
\end{remark}

\subsection{Relevant results}

The wave--Klein--Gordon system is a classical and important model in mathematical physics, and it has attracted considerable research interest in partial differential equations. Over the past few decades, there has been a vast literature on nonlinear wave equations, Klein--Gordon equations, and their coupled systems. The global existence and the long-time behavior of small-data solutions for the Cauchy problem are the focus of the present article.

In the three-dimensional setting, foundational techniques were first established for single wave or Klein--Gordon equations. In 1986, the global existence of solutions to 3D quasilinear wave equations with the null condition was proved independently by Christodoulou \cite{Christodoulou} and Klainerman \cite{Kla1984}, where Klainerman also introduced the vector field method \cite{Klainerman852}. Around the same time, for nonlinear Klein--Gordon equations in $(1+3)$-dimensional spacetime, Klainerman \cite{Klainerman85} proved small-amplitude global existence, while Shatah \cite{Shat85} developed the normal form method to handle quadratic nonlinearities; see also the work \cite{LindbladSoffer} by Lindblad and Soffer on asymptotic completeness. For coupled wave--Klein--Gordon systems, Georgiev \cite{Georgiev} first established global existence for small initial data in 1990. Later on, Ionescu and Pausader \cite{Ionescu-P} developed the Fourier based Z-norm method to prove global regularity for a coupled wave--Klein--Gordon system.

Due to the slow decay of linear waves, the two-dimensional case is substantially more delicate. In 2001, Alinhac \cite{Alinhac1} obtained the global existence for 2D quasilinear wave equations under the null condition. Subsequently, Dong, LeFloch, and Lei \cite{DoLeLe21} and Li \cite{Li21} independently established the uniform boundedness of the top-order energy for 2D quasilinear wave equations; see also the recent work \cite{HoTaYi} by Hou--Tao--Yin on quasilinear wave systems in the setting $\mathbb{R}^2 \times \mathbb{T}$. For the coupled two-dimensional wave--Klein--Gordon system, Dong \cite{Dong2005} proved global existence under the null condition in 2021. Recently, Dong--Ma--Yuan \cite{DongMaYuan} further studied the asymptotic behavior of the 2D coupled wave--Klein--Gordon system. In contrast to the three- and two-dimensional cases, the one-dimensional problem \eqref{eq:wave-KG} is solved in the present paper.

Last but not least, we briefly review some closely related results in geometric, gauge, and dispersive-wave coupled models where important advances have been achieved. LeFloch and Ma \cite{PLF-YM-cmp} used the hyperboloidal vector field method and established the global nonlinear stability of Minkowski space for the self-gravitating wave--Klein--Gordon model. Later on, Wang \cite{Wang} made significant progress on the Einstein--Klein--Gordon equations by developing an intrinsic hyperboloid approach, while Klainerman--Wang--Yang \cite{Klainerman-QW-SY} achieved notable progress for the massive Maxwell--Klein--Gordon equations by establishing global solutions. More recently, Ionescu and Pausader \cite{Ionescu-P2} established the global stability of Minkowski space for the Einstein--Klein--Gordon coupled system, and Dong--Guo--Li \cite{DongGuoLiKGZ} made important advances in the scattering theory for the two-dimensional Klein--Gordon--Zakharov system. For other one-dimensional models, we refer to Luli--Yang--Yu \cite{Luli} for semilinear wave systems with null conditions and to Ma \cite{Ma1D} for the cubic nonlinear wave--Klein--Gordon system.

\subsection{Main difficulties and proof strategy.}
There are several main difficulties in proving our result for the one-dimensional wave--Klein--Gordon system \eqref{eq:wave-KG}. First, the coupled wave--Klein--Gordon system lacks compatibility with the scaling vector field due to the mass term. 
We therefore use the Lorentz boost \(x\partial_t+t\partial_x\),
which commutes with both the wave and Klein--Gordon operators.
In hyperbolic coordinates this vector field is \(\partial_y\),
as shown in \eqref{eq:hyp-vector-fields}.
Second, as the spatial dimension decreases, the decay rates of solutions become slower, which makes lower-dimensional problems substantially more delicate and renders the one-dimensional case particularly difficult.

To overcome this, we first rewrite the wave equation in divergence form,
\[
  \partial_t(w_t+v w_x)=\partial_x(w_x+v w_t),
\]
which allows us to define a potential \(\phi\) by \(\phi_x=w_t+v w_x\) and \(\phi_t=w_x+v w_t\), with the constant fixed by setting \(\phi=0\) to the left of the initial support. We then introduce \(V=e^{\phi/2}v\) so that derivatives of the exponential cancel the quadratic coupling in the Klein--Gordon equation, leaving only terms of at least cubic order in \eqref{eq:gauged-kg}.

To extract the \(\tau^{-1/2}\) decay of \(V\), we set
\(\psi=\tau^{1/2}V\). This rescaling removes the first-order time
derivative from the linear part of the equation, since
\[
  V_{\tau\tau}+\tau^{-1}V_\tau
  =\tau^{-1/2}\left(\psi_{\tau\tau}
    +\frac{1}{4\tau^2}\psi\right).
\]
The source \(\mathcal R\) in the rescaled equation \eqref{eq:rescaled-kg}
satisfies \(\|\mathcal R\|_{H_y^{N-1}}\lesssim\varepsilon^3\tau^{-3/2}\)
under the low-order bootstrap bound, by Lemma~\ref{lem:kg-remainder}.
This bound is integrable in \(\tau\), so the low-order energy estimate
gives \(\|\psi\|_{H_y^1}\lesssim\varepsilon\) uniformly in time.
By one-dimensional Sobolev embedding and the bootstrap bound on
\(\|\phi\|_\infty\) in \eqref{eq:bootstrap-low},
the identity \(v=\tau^{-1/2}e^{-\phi/2}\psi\) then yields
\(\|v\|_\infty\lesssim\varepsilon\tau^{-1/2}\).

Another key issue is how to control the top-order wave derivatives while keeping the lower-order wave energy uniformly bounded. In the coordinates \((\tau,y)\), the vector \(U=(\phi_y,w_y)^{\mathsf{T}}\) satisfies \(\partial_\tau U=\tau^{-1}\partial_y(M(v)U)\) with the coefficient matrix \(M(v)\) in \eqref{eq:potential-system}. Lemma~\ref{lem:potential-system}
shows that \(H(v)=\operatorname{diag}(1,1-v^2)\) is uniformly positive
for \(|v|\le\tfrac12\) and that \(H(v)M(v)\) is symmetric. We therefore
use \(H(v)\) as the weight in the order-\(N\) wave energy. Integration by parts transfers the extra derivative in
\(\partial_y (\partial_y^{N}U)\) to \(HM\):
\[
  \frac{2}{\tau}\int_{\mathbb R} (\partial_y^N U)^{\mathsf T}HM\, (\partial_y^{N+1} U)\,dy
  =-\frac{1}{\tau}\int_{\mathbb R} (\partial_y^N U)^{\mathsf T}\partial_y(HM)(\partial_y^N U)\,dy .
\]
Thus, the leading transport term is controlled by the derivatives of \(U\) only up to order \(N\).

Fix \(0<\delta<\tfrac12\) and write \(A_m(\tau)=\|U(\tau,\cdot)\|_{H_y^m}\). At order \(N-1\), we use the Euclidean
energy \(A_{N-1}\). Its constant weight eliminates the term
involving \(H_\tau\). The skew-symmetric part of the transport
matrix now costs one derivative, contributing at most
\[
  C\tau^{-1}\|v\|_\infty^2 A_N(\tau)
  \lesssim \varepsilon^3\tau^{-3/2+\delta},
\]
where we have used \(\|v\|_\infty\lesssim\varepsilon\tau^{-1/2}\)
and the top-order bootstrap bound
\(A_N(\tau)\lesssim\varepsilon\tau^{1/2+\delta}\).
Since \(\delta<\tfrac12\), this contribution is integrable in \(\tau\).
See \eqref{eq:bootstrap-evol-low} for the differential inequality
and \eqref{eq:close-low-bound} for the resulting uniform bound.

Finally, we improve the bootstrap bound on \(\phi\) itself, since the wave energy only controls derivatives of \(\phi\). To achieve this, we integrate the equations for \(\phi\pm w\) along the characteristics of the system with coefficient matrix \(M(0)\). Since the source term is bounded in \(L^\infty_y\) by \(C\tau^{-1}\|v\|_\infty\|U\|_\infty\lesssim\varepsilon^2\tau^{-3/2}\) and is integrable, we obtain \(\|\phi(\tau)\|_\infty+\|w(\tau)\|_\infty\lesssim\varepsilon\). Together with the energy estimates, this closes the bootstrap and gives the hyperboloidal estimates for all \(\tau\ge2\).

To obtain the sharp decay rate, we use finite speed of propagation once more. For \(|y|>\log 2\), the ray with \(y\) fixed stays outside the support of \(\psi\) for \(2\le\tau<e^{|y|}\). Applying Duhamel's formula to the rescaled Klein--Gordon equation starting at \(\tau=e^{|y|}\) yields
\[
  |\psi(\tau,y)|\lesssim\varepsilon e^{-|y|/2}.
\]
Since \(\tau e^{|y|}=t+|x|\) and \(|x|\le t\) on the support, this gives the uniform decay estimate
\[
  |v(t,x)|=\tau^{-1/2}e^{-\phi/2}|\psi(\tau,y)|
    \lesssim \varepsilon (t+|x|)^{-1/2}
    \lesssim \varepsilon (1+t)^{-1/2}.
\]
See Corollary~\ref{cor:physical-time-decay} for details.

\subsection{Organization of the paper.}
Section~\ref{sec:pre} introduces the hyperbolic coordinates and records the energy inequalities used throughout the paper. Section~\ref{sec:hyperbolic-system} derives the first-order system for the potentials together with the rescaled Klein--Gordon equation. Section~\ref{sec:bootstrap} sets up the bootstrap assumptions and establishes the pointwise bounds. Section~\ref{sec:kg-sobolev} estimates the coefficients and source term of the rescaled equation in Sobolev spaces. Finally, Section~\ref{sec:bootstrap-closure} combines the wave and Klein--Gordon energy estimates with a characteristic bound for the potential to close the bootstrap and extend the solution globally.

\section{Preliminaries}\label{sec:pre}
\subsection{Hyperbolic coordinates}\label{sec:hyperbolic-coordinates}
We write \((t,x)\) for the Cartesian coordinates on \(\mathbb R^{1+1}\).
Our convention for the wave operator is
\[
  -\Box=\partial_t^2-\partial_x^2.
\]
Since the initial data at \(t=2\) are supported in \([-1,1]\),
finite propagation gives \(|x|\le t-1\) on the support of the solution.
This region lies in the forward cone \(t>|x|\). In this cone, we introduce the hyperbolic time \(\tau>0\)
and the rapidity \(y\in\mathbb R\) by
\begin{equation}
  t=\tau\cosh y,\qquad x=\tau\sinh y.
  \label{eq:hyp-coordinates}
\end{equation}
Since \(t\pm x=\tau e^{\pm y}\), the inverse transformation is
\begin{equation}
  \tau=\sqrt{t^2-x^2},\qquad
  y=\frac12\log\!\left(\frac{t+x}{t-x}\right).
  \label{eq:hyp-inverse}
\end{equation}

We use the same symbols \(w\) and \(v\) for the two fields in Cartesian coordinates \((t,x)\) and in hyperbolic coordinates \((\tau,y)\). By the chain rule, the scaling vector field and the Lorentz boost take the form
\begin{equation}
  \begin{aligned}
    \tau\partial_\tau
      &=\tau\cosh y\,\partial_t+\tau\sinh y\,\partial_x
       =t\partial_t+x\partial_x,\\
    \partial_y
      &=\tau\sinh y\,\partial_t+\tau\cosh y\,\partial_x
       =x\partial_t+t\partial_x.
  \end{aligned}
  \label{eq:hyp-vector-fields}
\end{equation}

Inverting these relations gives 
\[ \partial_t=\cosh y\,\partial_\tau-\frac{\sinh y}{\tau}\partial_y, \qquad \partial_x=-\sinh y\,\partial_\tau+\frac{\cosh y}{\tau}\partial_y. 
\]

Squaring the vector fields in \eqref{eq:hyp-vector-fields} and
subtracting gives
\[(\tau\partial_\tau)^2-\partial_y^2
  =\tau^2(\partial_t^2-\partial_x^2).
\]
Expanding \((\tau\partial_\tau)^2=\tau^2\partial_\tau^2+\tau\partial_\tau\), we obtain
\begin{equation}
  -\Box
  =\partial_\tau^2+\frac1\tau\partial_\tau-\frac1{\tau^2}\partial_y^2.
  \label{eq:hyp-wave-operator}
\end{equation}

For \(C^1\) scalar functions \(f\) and \(g\) of \((t,x)\),
we define the null forms \(Q_0\) and \(Q_{01}\) by
\[
  Q_0(f,g)=f_tg_t-f_xg_x,\qquad Q_{01}(f,g)=f_tg_x-f_xg_t.
\]
In the hyperbolic coordinates \((\tau,y)\), these null forms take the form
\begin{equation}\label{eq:hyp-null-forms}
\begin{aligned}
Q_0(f,g)&=f_\tau g_\tau-\tau^{-2}f_y g_y,\\
Q_{01}(f,g)&=\tau^{-1}\bigl(f_\tau g_y-f_y g_\tau\bigr).
\end{aligned}
\end{equation}

On each hyperbola we take spatial norms in \(y\) at fixed \(\tau\). We write \(\|f\|_{p}=\|f\|_{L^p(\mathbb R_y)}\).
For an integer \(m\ge0\), we use the norms
\[
  \|f\|_{H^m}^2=\sum_{j=0}^m\|\partial_y^j f\|_{2}^2,\qquad
  \|f\|_{W^{1,\infty}}=\|f\|_{\infty}+\|\partial_y f\|_{\infty} .
\]
For vector-valued functions we use the Euclidean norm, and for matrices
we use the induced operator norm. For a matrix or vector \(A\),
\(A^{\mathsf T}\) denotes its transpose. We write \(A\lesssim B\) if
\(A\le CB\) for a constant \(C\).

\subsection{Energy inequalities}

Throughout this subsection, we assume that the solutions and coefficients
are smooth, that the norms appearing below are finite and continuous in
\(\tau\), and that integration by parts in \(y\) produces no boundary terms.
We write
\[
  D_\tau^+ E(\tau)
  = \limsup_{h\to 0^+}\frac{E(\tau+h)-E(\tau)}{h}
\]
for the upper right Dini derivative. We also use a prime to denote differentiation with respect to \(\tau\).

\begin{lemma}[Energy inequality for the Klein--Gordon equation]
\label{lem:general-kg-energy}
Suppose \(u(\tau,y)\) satisfies the equation
\[
  u_{\tau\tau}-p(\tau)u_{yy}+q(\tau)u=f,
\]
where \(p,q>0\) and \(p',q'\le0\). For an integer \(m\ge0\), define the energy
\[
  E_m^2=\|u_\tau\|_{H^m}^2+p(\tau)\|u_y\|_{H^m}^2+q(\tau)\|u\|_{H^m}^2.
\]
Then
\begin{equation}
  D_\tau^+E_m\le\|f\|_{H^m}.
  \label{eq:general-kg-energy}
\end{equation}
The estimate also holds in integral form.
\end{lemma}

\begin{proof}
Fix \(0\le j\le m\) and set \(z=\partial_y^j u\). Since \(p\) and \(q\)
depend only on \(\tau\), applying \(\partial_y^j\) to the equation gives
\[
  z_{\tau\tau}-p z_{yy}+qz=\partial_y^j f.
\]
Multiplying both sides by \(2z_\tau\) and integrating over \(\mathbb{R}_y\) give
\[
  \frac{d}{d\tau}\|z_\tau\|_{2}^2
  +p\frac{d}{d\tau}\|z_y\|_{2}^2
  +q\frac{d}{d\tau}\|z\|_2^2
  =2\int_{\mathbb{R}}z_\tau\,\partial_y^j f\,dy,
\]
in which we used
\[-2p\int_{\mathbb{R}} z_\tau z_{yy}\,dy
  =2p\int_{\mathbb{R}} z_{\tau y}z_y\,dy.\]
Using the product rule and summing over \(j=0,\ldots,m\), we obtain
\[
  (E_m^2)'=p'\|u_y\|_{H^m}^2+q'\|u\|_{H^m}^2
           +2\sum_{j=0}^m\int_{\mathbb{R}}
             (\partial_\tau\partial_y^j u)(\partial_y^j f)\,dy.
\]
Applying Cauchy--Schwarz in \(y\) and then in \(j\) bounds the last term by
\[2\sum_{j=0}^m\int_{\mathbb{R}}
             (\partial_\tau\partial_y^j u)(\partial_y^j f)\,dy\le2\|u_\tau\|_{H^m}\|f\|_{H^m}.\]
Since \(p',q'\le 0\) and
\(\|u_\tau\|_{H^m}\le E_m\), it follows that \((E_m^2)'\le 2E_m\|f\|_{H^m}\).
For \(\rho>0\), this implies
\[
  \frac{d}{d\tau}(E_m^2+\rho^2)^{1/2}\le\|f\|_{H^m}.
\]
Integrating and letting \(\rho\to 0^+\) yields the integral estimate, and taking right difference quotients gives \eqref{eq:general-kg-energy}.
\end{proof}

\begin{lemma}[Energy inequality for first-order systems]
\label{lem:general-system-energy}
Suppose that \(W(\tau,y)\) satisfies
\[
  W_\tau=P(\tau,y)W_y+F.
\]
Let \(G=G^{\mathsf T}\) be a matrix satisfying
\(\lambda I\le G(\tau,y)\le\Lambda I\) for constants \(0<\lambda\le\Lambda\).
Write \(GP=S+J\), where \(S\) and \(J\) are its symmetric and
skew-symmetric parts, respectively. Setting
\[
  E_G^2=\int_{\mathbb{R}}W^{\mathsf T}GW\,dy,
\]
we have
\begin{equation}
  D_\tau^+E_G
  \le \frac{\|G_\tau-S_y\|_\infty}{2\lambda}\,E_G
      +\frac{\|J\|_\infty}{\sqrt{\lambda}}\,\|W_y\|_2
      +\sqrt{\Lambda}\,\|F\|_2.
  \label{eq:general-system-energy}
\end{equation}
The matrix norms are induced by the Euclidean norm. The estimate also
holds in integral form.
\end{lemma}

\begin{proof}
Multiply \(W_\tau=P W_y+F\) on the left by \(2W^{\mathsf T}G\).
Using the product rule, the symmetry of \(G\), and the decomposition
\(GP=S+J\), we obtain
\[
  \partial_\tau(W^{\mathsf T}GW)-W^{\mathsf T}G_\tau W
  =2W^{\mathsf T}S W_y+2W^{\mathsf T}J W_y+2W^{\mathsf T}GF .
\]
Since \(S=S^{\mathsf T}\), the first term on the right is
\(\partial_y(W^{\mathsf T}SW)-W^{\mathsf T}S_yW\).
Integrating over \(\mathbb R_y\) eliminates the spatial divergence and yields
\[
  (E_G^2)'=\int_{\mathbb R}\!\left[
  W^{\mathsf T}(G_\tau-S_y)W+2W^{\mathsf T}J W_y
  +2W^{\mathsf T}GF\right]dy .
\]
For the first two terms we use \(\|W\|_2\le E_G/\sqrt{\lambda}\),
which follows from \(G\ge\lambda I\). For the source term, the positive
symmetric square root of \(G\) satisfies
\[
  \|G^{1/2}F\|_2^2=\int_{\mathbb R}F^{\mathsf T}GF\,dy
  \le \Lambda\|F\|_2^2,
\]
since \(G\le\Lambda I\). Cauchy--Schwarz applied to \(G^{1/2}W\) and
\(G^{1/2}F\), together with \(\|G^{1/2}W\|_2=E_G\), thus gives
\[
  \left|2\int_{\mathbb R}W^{\mathsf T}GF\,dy\right|
  \le 2E_G\sqrt{\Lambda}\,\|F\|_2 .
\]
Combining the three bounds yields
\[
  (E_G^2)'\le 2E_G\left(
    \frac{\|G_\tau-S_y\|_\infty}{2\lambda}\,E_G
    +\frac{\|J\|_\infty}{\sqrt{\lambda}}\,\|W_y\|_2
    +\sqrt{\Lambda}\,\|F\|_2\right).
\]
For \(\rho>0\), divide by \(2\sqrt{E_G^2+\rho^2}\) and use
\(E_G/\sqrt{E_G^2+\rho^2}\le1\). Integrating and letting
\(\rho\to 0^+\) gives the integral estimate. Taking right
difference quotients yields \eqref{eq:general-system-energy}.
\end{proof}

\section{The coupled system in hyperbolic coordinates}\label{sec:hyperbolic-system}


We reformulate the coupled system to obtain the evolution equations
used in the wave and Klein--Gordon energy estimates. In Cartesian
coordinates \((t,x)\), the wave equation
\(w_{tt}-w_{xx}=w_tv_x-w_xv_t\) has the conservation form
\[
  \partial_t(w_t+v w_x)=\partial_x(w_x+v w_t).
\]
This identity is the compatibility condition for the existence of a
potential \(\phi\) satisfying
\[
  \phi_x=w_t+v w_x,\qquad \phi_t=w_x+v w_t.
\]
We fix its additive constant by taking \(\phi=0\) to the left of the
initial support on \(t=2\).

We use the potential identities to derive a first-order system for \((\phi,w)\) in the coordinates \((\tau,y)\). We then use \(\phi\) to eliminate the quadratic coupling in the Klein--Gordon equation.

\subsection{The first-order system for the potentials}
\label{sec:potential-system}
We group \(w\) with its potential by setting \(X=(\phi,w)^{\mathsf T}\) and
\(U=X_y=(U_1,U_2)^{\mathsf T}\). The wave energy is expressed in terms of \(U\), allowing \(\phi\) to approach different constants as \(y\to\pm\infty\). 
By \eqref{eq:hyp-vector-fields}, the potential relations become
\begin{equation}
  \phi_y=\tau w_\tau+vw_y,\qquad \tau\phi_\tau=w_y+\tau vw_\tau.
  \label{eq:potential-identities}
\end{equation}

\begin{lemma}[First-order system and symmetrizer]
\label{lem:potential-system}
Under \eqref{eq:potential-identities}, the potentials and their
spatial derivatives satisfy
\begin{equation}
  X_\tau=\frac{1}{\tau}M(v)X_y,\qquad U_\tau=\frac{1}{\tau}\partial_y\bigl(M(v)U\bigr),
  \qquad
  M(v)=\begin{pmatrix}v&1-v^2\\1&-v\end{pmatrix}.
  \label{eq:potential-system}
\end{equation}
With \(H(v)=\operatorname{diag}(1,1-v^2)\), the coefficient matrix satisfies
\begin{align}
  M(v)^2 &= I_2, \label{eq:potential-involution}\\
  H(v)M(v) &= \bigl(H(v)M(v)\bigr)^{\mathsf T},
    \label{eq:potential-symmetrizer}\\
  M(v)-M(v)^{\mathsf T}
    &= v^2\begin{pmatrix}0&-1\\1&0\end{pmatrix}.
    \label{eq:potential-skew}
\end{align}
In particular, \(H(v)\) is a uniformly positive symmetrizer for
\(|v|\le\tfrac12\).
\end{lemma}

\begin{proof}
$\bullet$ \emph{The first-order system.} Combining 
\(\tau\phi_\tau=w_y+v\tau w_\tau\) with the definition of 
\(U\), we obtain 
\[ \tau\phi_\tau=vU_1+(1-v^2)U_2,\qquad 
   \tau w_\tau=U_1-vU_2. 
\] 
These are the two rows of \(\tau X_\tau=M(v)U\). Differentiating in \(y\) gives the equation for \(U\).

$\bullet$ \emph{Matrix identities.} Direct multiplication gives
\[
  H(v)M(v)=
  \begin{pmatrix}
    v&1-v^2\\
    1-v^2&-v(1-v^2)
  \end{pmatrix},
\]
which is symmetric. Squaring \(M(v)\) gives \(I_2\), and subtracting
its transpose gives \eqref{eq:potential-skew}. 

$\bullet$ \emph{Positivity of \(H(v)\).} For
\(|v|\le\tfrac12\) and \(z\in\mathbb R^2\),
\[
  \tfrac34|z|^2\le z^{\mathsf T}H(v)z
    =z_1^2+(1-v^2)z_2^2\le |z|^2. \qedhere
\]
\end{proof}

\subsection{The rescaled Klein--Gordon equation}
We first introduce \(V=e^{\phi/2}v\). By the potential identities,
the quadratic term cancels in the equation for \(V\),
leaving only cubic and higher-order terms on the right-hand side.
We then set \(\psi=\tau^{1/2}V\) to remove the first-order \(\tau\)-derivative
from the linear part of the equation. 
The coefficients of the resulting remainder \(\mathcal R\) in the rescaled equation depend on \(v\) and the wave derivatives
\(\tau w_\tau\) and \(w_y\).

\begin{lemma}[The gauged and rescaled Klein--Gordon equations]\label{lem:hyperbolic-kg}
For a smooth solution \((w,v)\) of \eqref{eq:wave-KG} with associated
potential \(\phi\), set
\[
  V=e^{\phi/2}v,\qquad \psi=\tau^{1/2}V,
\]
so that \(v=\tau^{-1/2}e^{-\phi/2}\psi\). Then \(V\) satisfies
\begin{equation}
  (-\Box+1)V
  =\frac32 v Q_0(w,V)+\frac12 v^2 Q_{01}(w,V)
   -\frac{1+3v^2}{4}\,V Q_0(w,w)
  \label{eq:gauged-kg}
\end{equation}
and, under the rescaling \(\psi=\tau^{1/2}V\), this equation is equivalent to
\begin{equation}
  \psi_{\tau\tau}-\frac{1}{\tau^{2}}\psi_{yy}
  +\left(1+\frac{1}{4\tau^2}\right)\psi=\mathcal R,\qquad
  \mathcal R=\frac{\alpha}{\tau}\psi_\tau+\frac{\beta}{\tau^2}\psi_y
             -\frac{\gamma}{\tau^2}\psi,
  \label{eq:rescaled-kg}
\end{equation}
where
\begin{equation}
  \begin{aligned}
    \alpha &= \frac32 \tau vw_\tau-\frac12 v^2w_y, &
    \beta  &= -\frac32 vw_y+\frac12 \tau v^2w_\tau,\\
    \gamma &= \frac{\alpha}{2}+\frac{1+3v^2}{4}(\tau^2w_\tau^2-w_y^2).
  \end{aligned}
  \label{eq:rescaled-kg-coefficients}
\end{equation}
\end{lemma}

\begin{proof}
$\bullet$ \emph{The equation for \(V\).}
Applying the product rule to \(V=e^{\phi/2}v\) gives
\[
  (-\Box+1)V=e^{\phi/2}\left[
    (-\Box+1)v+Q_0(\phi,v)
    +\frac v2(-\Box\phi)+\frac v4 Q_0(\phi,\phi)\right].
\]
For the first two terms, substituting the potential identity
\[
  Q_0(\phi,v)=-Q_{01}(w,v)+vQ_0(w,v)
\]
and the Klein--Gordon equation
\[
  (-\Box+1)v=Q_{01}(w,v)
\]
cancels \(Q_{01}(w,v)\) and leaves the cubic term \(vQ_0(w,v)\).

Differentiating the potential identities
and using the wave equation gives
\[
  -\Box\phi
  =\partial_t(w_x+vw_t)-\partial_x(w_t+vw_x)
  =Q_0(w,v)+vQ_{01}(w,v).
\]
Since \(Q_0(\phi,\phi)=(v^2-1)Q_0(w,w)\), we obtain
\[
  (-\Box+1)V
  =e^{\phi/2}\left[\frac32 v Q_0(w,v)+\frac12 v^2 Q_{01}(w,v)
    +\frac{v(v^2-1)}4 Q_0(w,w)\right].
\]
To express the differentiated factors in terms of \(V\), use
\(Q_0(w,\phi)=vQ_0(w,w)\) and \(Q_{01}(w,\phi)=Q_0(w,w)\). Then
\[
  \begin{aligned}
    e^{\phi/2}Q_0(w,v)&=Q_0(w,V)-\tfrac12 vV Q_0(w,w),\\
    e^{\phi/2}Q_{01}(w,v)&=Q_{01}(w,V)-\tfrac12 V Q_0(w,w).
  \end{aligned}
\]
Substituting these identities and the expressions for \(-\Box\phi\) and
\(Q_0(\phi,\phi)\) into the formula for \((-\Box+1)V\) yields the
coefficient \(-(1+3v^2)/4\) of \(VQ_0(w,w)\) and proves
\eqref{eq:gauged-kg}.

$\bullet$ \emph{The equation for \(\psi\).}
Substituting \(V=\tau^{-1/2}\psi\) into \eqref{eq:hyp-wave-operator}, we obtain
\begin{align*}
  \tau^{1/2}(-\Box+1)V
  &=\tau^{1/2}(\partial_\tau^2+\tau^{-1}\partial_\tau-\tau^{-2}\partial_y^2)(\tau^{-1/2}\psi)+\psi\\
  &=\psi_{\tau\tau}-\frac{1}{\tau^{2}}\psi_{yy}+\left(1+\frac{1}{4\tau^2}\right)\psi.
\end{align*}
On the right-hand side of \eqref{eq:gauged-kg}, the vector fields in
\eqref{eq:hyp-vector-fields} give
\[
\begin{aligned}
  \tau^2Q_0(w,V)&=\tau^2w_\tau V_\tau-w_yV_y,\\
  \tau^2Q_{01}(w,V)&=\tau w_\tau V_y-\tau w_yV_\tau,\\
  \tau^2Q_0(w,w)&=\tau^2w_\tau^2-w_y^2.
\end{aligned}
\]
Rewriting the right-hand side of \eqref{eq:gauged-kg} via these identities,
we see that the coefficients of \(V_\tau\) and \(V_y\) are
\(\tau^{-1}\alpha\) and \(\tau^{-2}\beta\), respectively. Multiplying \eqref{eq:gauged-kg} by \(\tau^{1/2}\) and substituting
\[
  \tau^{1/2}V_\tau=\psi_\tau-\tfrac1{2\tau}\psi
  \qquad\text{and}\qquad
  \tau^{1/2}V_y=\psi_y
\]
yields \eqref{eq:rescaled-kg}. Since \(\tau>0\), the substitution is invertible, which proves the equivalence.
\end{proof}

\section{Bootstrap estimates}\label{sec:bootstrap}
The wave energy estimates require pointwise control of \(v_y\) and
\(vv_\tau\), which occur when we differentiate the transport matrix
and its symmetrizer. We obtain these bounds from a uniform smallness
assumption one derivative below the top order. At the top order we
allow the wave and Klein--Gordon energies to grow at different rates.

\subsection{Energy norms}

Fix the top order \(N\ge6\). For each integer \(m\ge0\), let \(A_m\) measure the
wave variable \(U=X_y\) and let \(K_m\) be the energy of the
rescaled Klein--Gordon field \(\psi\) on the hyperbola of radius \(\tau\).
Explicitly,
\[
  \begin{aligned}
    A_m(\tau)&=\|U(\tau)\|_{H^m},\\
    K_m(\tau)^2&=\|\psi_\tau(\tau)\|_{H^m}^2+\frac{1}{\tau^{2}}\|\psi_y(\tau)\|_{H^m}^2
       +\left(1+\frac{1}{4\tau^2}\right)\|\psi(\tau)\|_{H^m}^2 .
  \end{aligned}
\]
At order \(N\), the symmetrizer \(H(v)\) from
Lemma~\ref{lem:potential-system} gives the wave energy
\[
  B_N(\tau)^2=\sum_{j=0}^N\int_{\mathbb R}
     (\partial_y^jU)^{\mathsf T}H(v)\,\partial_y^jU\,dy .
\]
The energies \(B_N(\tau)\) and \(A_N(\tau)\) are uniformly equivalent whenever
\(\|v(\tau)\|_\infty\le\tfrac12\).

\subsection{Bootstrap assumptions and pointwise bounds}
\label{sec:bootstrap-pointwise}
For the Cauchy problem in Theorem~\ref{thm:main-result}, take
\(\tau_0=2\). For sufficiently small \(\varepsilon\), standard local
well-posedness gives a smooth solution on \(2\le t\le3\).
By finite propagation, its support lies in \(|x|\le t-1\). On the
portion of the initial hyperbola \(t^2-x^2=4\) in this region,
\[
  4=t^2-x^2\ge t^2-(t-1)^2=2t-1,
  \qquad 2\le t\le\tfrac52 .
\]
Thus the local solution covers this entire arc. Moreover,
\(2e^{-|y|}=t-|x|\ge1\) there, so \(|y|\le\log 2\).

On this compact arc, the local energy estimates following from
\eqref{eq:initial-data-smallness} and one-dimensional Sobolev
embedding control the derivatives of \(v\) and \(w\) through order
\(N+1\) by \(C_{N}\varepsilon\). The potential identities give
the corresponding bounds for the derivatives of \(\phi\), and its
normalization on the left gives the same bound for \(\phi\) itself.
Outside this arc, \(v\) and \(w\) vanish and \(\phi\) is constant
on each exterior component. Thus, by the definitions of \(A_N\),
\(K_N\), and \(X\),
\begin{equation}
  A_N(\tau_0)+K_N(\tau_0)+\|X(\tau_0)\|_\infty
  \le C_0\varepsilon,
  \label{eq:bootstrap-data}
\end{equation}
where \(C_0\) depends only on \(N\).


Fix \(0<\delta<\tfrac12\).
Choose a bootstrap constant \(C_1\ge1\) sufficiently large relative to \(C_0\)
and take \(\varepsilon>0\) so small that \(C_1\varepsilon\le1\).
Let \([\tau_0,\tau^*]\) be an interval of existence and suppose that, for every
\(\tau\in[\tau_0,\tau^*]\),
\begin{align}
  A_{N-1}(\tau)+K_{N-1}(\tau)+\|X(\tau)\|_{\infty}
    &\le C_1\varepsilon, \label{eq:bootstrap-low}\\
  B_N(\tau)&\le C_1\varepsilon \tau^{1/2+\delta}, \label{eq:bootstrap-wave}\\
  K_N(\tau)&\le C_1\varepsilon \tau^\delta. \label{eq:bootstrap-kg}
\end{align}

\begin{proposition}[Pointwise estimates from the low-order bootstrap bound]
\label{prop:bootstrap-pointwise}
For every \(\tau\in[\tau_0,\tau^*]\), the low-order bound \eqref{eq:bootstrap-low}
implies
\begin{align}
  \|v(\tau)\|_{W^{1,\infty}}&\lesssim C_1\varepsilon \tau^{-1/2},
    \label{eq:pointwise-v}\\
  \|U(\tau)\|_{W^{1,\infty}}&\lesssim C_1\varepsilon,
    \label{eq:pointwise-U}\\
  \|v_\tau(\tau)\|_{\infty}&\lesssim C_1\varepsilon \tau^{-1/2},
    \label{eq:pointwise-v-tau}\\
  \|v(\tau)v_\tau(\tau)\|_{\infty}&\lesssim (C_1\varepsilon)^2\tau^{-1}.
    \label{eq:pointwise-v-v-tau}
\end{align}
The implicit constants are independent of \(C_1\), \(\varepsilon\), and \(\tau^*\).
\end{proposition}

\begin{proof}
 
$\bullet$ \emph{Bound for \(U\).} By the definitions of \(A_{N-1}\) and
\(K_{N-1}\), the low-order bound and the one-dimensional Sobolev
embeddings \(H^2\hookrightarrow W^{1,\infty}\) and
\(H^1\hookrightarrow L^\infty\) give
\[
  \|U\|_{W^{1,\infty}}+\|\psi\|_{W^{1,\infty}}
  +\|\psi_\tau\|_{\infty}\lesssim C_1\varepsilon.
\]
This proves \eqref{eq:pointwise-U}. 

$\bullet$ \emph{Bound for \(v\).} Since \(\|\phi\|_\infty\le C_1\varepsilon\le1\)
and \(\phi_y=U_1\), the identities
\[
  v=\tau^{-1/2}e^{-\phi/2}\psi,\qquad
  v_y=\tau^{-1/2}e^{-\phi/2}
       \left(\psi_y-\tfrac12 U_1\psi\right)
\]
give \eqref{eq:pointwise-v}.

$\bullet$ \emph{Bound for \(v_\tau\).} The first row of \eqref{eq:potential-system} yields
\(\tau\phi_\tau=vU_1+(1-v^2)U_2\). Hence
\(\|\phi_\tau\|_\infty\lesssim C_1\varepsilon\tau^{-1}\) by
\eqref{eq:pointwise-v}--\eqref{eq:pointwise-U}.
Differentiating \(v=\tau^{-1/2}e^{-\phi/2}\psi\) in \(\tau\) gives
\[
  v_\tau=\tau^{-1/2}e^{-\phi/2}
    \left(\psi_\tau-\tfrac12(\tau^{-1}+\phi_\tau)\psi\right).
\]
Consequently \(\|v_\tau\|_\infty\lesssim
C_1\varepsilon(\tau^{-1/2}+\tau^{-3/2})\lesssim C_1\varepsilon \tau^{-1/2}\), which proves
\eqref{eq:pointwise-v-tau}. 

$\bullet$ \emph{Bound for \(vv_\tau\).} Combining \eqref{eq:pointwise-v} and \eqref{eq:pointwise-v-tau} gives
\[
  \|v v_\tau\|_\infty
  \le \|v\|_\infty \|v_\tau\|_\infty
  \lesssim (C_1\varepsilon)^2\tau^{-1}.
\]
\end{proof}

\section{Nonlinear estimates for the Klein--Gordon equation}\label{sec:kg-sobolev}
We use Sobolev bounds for \(v\) and the wave energy \(A_m\) to estimate
the coefficients of the rescaled Klein--Gordon equation and obtain the
factor \(\tau^{-3/2}\) in the bound for its remainder \(\mathcal R\).
Throughout this section we assume the low-order bootstrap bound
\eqref{eq:bootstrap-low} and work at fixed \(\tau\in[\tau_0,\tau^*]\).

\subsection{Sobolev estimates for the field \(v\)}

We recover \(v=\tau^{-1/2}e^{-\phi/2}\psi\) from the rescaled field \(\psi\).
The coefficients of the rescaled equation require bounds for \(v\)
through order \(N\). The top-order wave estimate also requires a bound
at order \(N+1\).

\begin{lemma}
\label{lem:v-sobolev}
For \(m\in\{N-1,N\}\),
\begin{equation}
  \|v\|_{H^m}\lesssim \tau^{-1/2}
     \bigl(K_m+C_1\varepsilon A_{m-1}\bigr).
  \label{eq:v-sobolev}
\end{equation}
\end{lemma}

\begin{proof}
Write \(v=\tau^{-1/2}g\psi\), where \(g=e^{-\phi/2}\).
By \eqref{eq:bootstrap-low}, \(\|\phi\|_\infty\le C_1\varepsilon\le1\),
so \(\|g\|_\infty\le e^{1/2}\). The same bound and the embedding
\(H^1(\mathbb R)\hookrightarrow L^\infty(\mathbb R)\) give
\(\|\psi\|_\infty\lesssim K_{N-1}\le C_1\varepsilon\).
For \(1\le k\le N+1\), the tame product estimate gives
\[
  \|g\psi\|_{H^k}\lesssim
  \|g\|_\infty\|\psi\|_{H^k}
  +\|\psi\|_\infty\|\partial_y g\|_{H^{k-1}} .
\]
Here the undifferentiated factor \(g\) is only measured in \(L^\infty\),
which allows \(\phi\) to have different limits at the two ends of the line.

To estimate \(\partial_y g\), observe that each derivative
\(\partial_y^j g\), \(1\le j\le k\), is a sum of terms of the form
\(C g\prod_{i=1}^r\partial_y^{j_i}\phi\), with
\(j_i\ge1\) and \(\sum_i j_i=j\).
Choose \(i_*\) with \(j_{i_*}=\max_i j_i\). By H\"older's inequality,
we estimate \(\partial_y^{j_{i_*}}\phi\) in \(L^2\) and all other
factors, including \(g\), in \(L^\infty\). Since \(\phi_y=U_1\),
\(\|\partial_y^{j_{i_*}}\phi\|_2\le\|U_1\|_{H^{k-1}}\).
Each remaining derivative has order
\(\ell\le\lfloor j/2\rfloor\le N-1\), so
\(\|\partial_y^\ell\phi\|_\infty
\lesssim\|\partial_y^{\ell-1}U_1\|_{H^1}\le A_{N-1}\le1\).
Thus \(\|\partial_y g\|_{H^{k-1}}\lesssim A_{k-1}\), and we obtain
\begin{equation}
  \|v\|_{H^k}\lesssim \tau^{-1/2}
    \bigl(\|\psi\|_{H^k}+C_1\varepsilon A_{k-1}\bigr),
  \qquad 1\le k\le N+1 .
  \label{eq:v-multiplier}
\end{equation}
For \(k=m\), the mass term in \(K_m^2\) gives
\(\|\psi\|_{H^m}\le K_m\), proving \eqref{eq:v-sobolev}.
\end{proof}

\begin{lemma}
\label{lem:v-sobolev-top}
At order \(N+1\),
\begin{equation}
  \|v\|_{H^{N+1}}\lesssim \tau^{1/2}K_N
     +C_1\varepsilon \tau^{-1/2}A_N .
  \label{eq:v-sobolev-top}
\end{equation}
\end{lemma}

\begin{proof}
The spatial derivative term in \(K_N^2\) carries the weight \(\tau^{-2}\).
Since \(\tau\ge1\), the energy definition gives
\[
  \tau^{-2}\|\psi\|_{H^{N+1}}^2
  =\tau^{-2}\|\psi\|_2^2+\tau^{-2}\|\psi_y\|_{H^N}^2
  \le K_N^2 .
\]
Taking \(k=N+1\) in \eqref{eq:v-multiplier} and inserting
\(\|\psi\|_{H^{N+1}}\le \tau K_N\) yields \eqref{eq:v-sobolev-top}.
\end{proof}

\subsection{Coefficient and remainder estimates for the rescaled equation}

\begin{lemma}[Pointwise coefficient bounds]
\label{lem:coefficients-pointwise}
The coefficients in \eqref{eq:rescaled-kg-coefficients} satisfy
\begin{equation}
  \|\alpha\|_{\infty}+\|\beta\|_{\infty}
    \lesssim (C_1\varepsilon)^2\tau^{-1/2},\qquad
  \|\gamma\|_{\infty}\lesssim (C_1\varepsilon)^2.
  \label{eq:coefficients-pointwise}
\end{equation}
\end{lemma}

\begin{proof}
By
Proposition~\ref{prop:bootstrap-pointwise},
\(\|v\|_\infty\lesssim C_1\varepsilon \tau^{-1/2}\) and \(\|U\|_\infty\lesssim C_1\varepsilon\).
Since \(\tau w_\tau=U_1-vU_2\) and \(w_y=U_2\), it follows that
\(\|\tau w_\tau\|_\infty+\|w_y\|_\infty\lesssim C_1\varepsilon\).
Each term in \(\alpha\) and \(\beta\) contains a factor \(v\) or \(v^2\)
multiplying \(\tau w_\tau\) or \(w_y\). This gives the first estimate since \(\tau\ge1\).

For the second, use
\(\gamma=\alpha/2+(1+3v^2)(\tau^2w_\tau^2-w_y^2)/4\) and
\(\|\tau^2w_\tau^2-w_y^2\|_\infty\lesssim(C_1\varepsilon)^2\).
\end{proof}

\begin{lemma}[Sobolev coefficient bounds]
\label{lem:coefficients-sobolev}
For \(m\in\{N-1,N\}\),
\begin{align}
  \|\alpha\|_{H^m}+\|\beta\|_{H^m}
    &\lesssim C_1\varepsilon \tau^{-1/2}(K_m+A_m),
      \label{eq:alpha-beta-sobolev}\\
  \|\gamma\|_{H^m}
    &\lesssim C_1\varepsilon A_m+C_1\varepsilon \tau^{-1/2}(K_m+A_m).
      \label{eq:gamma-sobolev}
\end{align}
\end{lemma}

\begin{proof}
Lemma~\ref{lem:v-sobolev} and
\(A_{m-1}\le A_m\) give \(\|v\|_{H^m}\lesssim \tau^{-1/2}(K_m+A_m)\).
Applying the tame product estimate to \(\tau w_\tau=U_1-vU_2\) and using \(w_y=U_2\)
gives
\begin{align*}
    \|\tau w_\tau\|_{H^m}+\|w_y\|_{H^m}
  &\lesssim (1+\|v\|_\infty)A_m+\|U_2\|_\infty\|v\|_{H^m}\\
  &\lesssim A_m+C_1\varepsilon\|v\|_{H^m}\\
  &\lesssim A_m+C_1\varepsilon \tau^{-1/2}(K_m+A_m).
\end{align*}
  
$\bullet$ \emph{Bounds for \(\alpha\) and \(\beta\).} For \(h=\tau w_\tau,w_y\) and \(p=1,2\), the same product estimate yields
\[
  \|v^p h\|_{H^m}
  \lesssim \|v\|_\infty^p\|h\|_{H^m}
       +\|h\|_\infty\|v\|_\infty^{p-1}\|v\|_{H^m}
  \lesssim C_1\varepsilon \tau^{-1/2}(K_m+A_m).
\]
This proves \eqref{eq:alpha-beta-sobolev} by the formulas for
\(\alpha\) and \(\beta\).

$\bullet$ \emph{Bound for \(\gamma\).} Set \(D=\tau^2w_\tau^2-w_y^2\). The preceding bounds imply
\[
  \|D\|_\infty\lesssim (C_1\varepsilon)^2,\qquad
  \|D\|_{H^m}\lesssim C_1\varepsilon(\|\tau w_\tau\|_{H^m}+\|w_y\|_{H^m})
  \lesssim C_1\varepsilon A_m+(C_1\varepsilon)^2\tau^{-1/2}(K_m+A_m).
\]
The term \(D/4\) in \(\gamma\) thus contributes \(C_1\varepsilon A_m\) even
without a factor of \(v\). For the full quadratic contribution we have
\begin{align*}
    \|D+3v^2D\|_{H^m}
  &\lesssim (1+\|v\|_\infty^2)\|D\|_{H^m}
      +\|D\|_\infty\|v\|_\infty\|v\|_{H^m}\\
  &\lesssim C_1\varepsilon A_m+(C_1\varepsilon)^2\tau^{-1/2}(K_m+A_m).
\end{align*}

Adding the bound for \(\alpha/2\) and using \(C_1\varepsilon\le1\) gives
\eqref{eq:gamma-sobolev}.
\end{proof}

\begin{lemma}[Estimate for the Klein--Gordon remainder]
\label{lem:kg-remainder}
For \(m\in\{N-1,N\}\), the remainder \(\mathcal R\) in
\eqref{eq:rescaled-kg} satisfies
\[
  \|\mathcal R\|_{H^m}
  \lesssim (C_1\varepsilon)^2\tau^{-3/2}(K_m+A_m).
\]
\end{lemma}

\begin{proof}
Set \(D_m=K_m+A_m\). The definition of
\(K_m\) gives
\[
  \|\psi_\tau\|_{H^m}\le K_m,\qquad
  \|\psi_y\|_{H^m}\le \tau K_m,\qquad
  \|\psi\|_{H^m}\le K_m.
\]
At low order, \eqref{eq:bootstrap-low} and Sobolev embedding give
\(\|\psi_\tau\|_\infty+\|\psi_y\|_\infty+\|\psi\|_\infty
\lesssim K_{N-1}\le C_1\varepsilon\).
Applying the tame product estimate to the first two terms of
\(\mathcal R\) and using Lemmas~\ref{lem:coefficients-pointwise}
and~\ref{lem:coefficients-sobolev}, we obtain
\begin{align*}
  \frac1\tau\|\alpha\psi_\tau\|_{H^m}
    &\lesssim \frac1\tau\bigl((C_1\varepsilon)^2\tau^{-1/2}K_m
              +(C_1\varepsilon)^2 \tau^{-1/2}D_m\bigr)
     \lesssim (C_1\varepsilon)^2\tau^{-3/2}D_m,\\
  \frac1{\tau^2}\|\beta\psi_y\|_{H^m}
    &\lesssim \frac1{\tau^2}\bigl((C_1\varepsilon)^2\tau^{-1/2}(\tau K_m)
              +(C_1\varepsilon)^2 \tau^{-1/2}D_m\bigr)
     \lesssim (C_1\varepsilon)^2\tau^{-3/2}D_m .
\end{align*}
For the last term, \eqref{eq:gamma-sobolev} and \(\tau\ge1\) give
\(\|\gamma\|_{H^m}\lesssim C_1\varepsilon D_m\). Together with the pointwise
bound \(\|\gamma\|_\infty\lesssim (C_1\varepsilon)^2\), this yields
\[
  \frac1{\tau^2}\|\gamma\psi\|_{H^m}
  \lesssim \tau^{-2}\bigl((C_1\varepsilon)^2 K_m+(C_1\varepsilon)^2 D_m\bigr)
  \lesssim (C_1\varepsilon)^2 \tau^{-2}D_m .
\]
Summing the three bounds and using \(\tau^{-2}\le \tau^{-3/2}\)
proves the claim.
\end{proof}

\section{Bootstrap improvement and global existence}\label{sec:bootstrap-closure}
We combine the pointwise estimates of Proposition~\ref{prop:bootstrap-pointwise}
with the bound for \(\mathcal R\) in Lemma~\ref{lem:kg-remainder} to
derive energy inequalities on \([\tau_0,\tau^*]\). In the lower-order wave
estimate, the factor \(\tau^{-2}\) makes the contribution of the highest
wave derivatives integrable for \(\delta<1/2\). At the top wave
order, a weighted Gronwall estimate controls the growth due to the
coefficient of size \((C_1\varepsilon)^2/\tau\). Together with the characteristic estimate for \(X\), these inequalities improve the three bootstrap bounds uniformly
in \(\tau^*\). The local continuation criterion then gives global existence. We then use the support of the solution and the low-order bounds to
derive the \(L^\infty_x\) decay of the Klein--Gordon field on each time slice
\(t\ge2\).

Choosing \(\varepsilon\) sufficiently small so that
\eqref{eq:pointwise-v} yields \(\|v\|_\infty \le 1/2\),
Lemma~\ref{lem:potential-system} gives
\begin{equation}
  \tfrac{\sqrt{3}}{2} A_N \le B_N \le A_N.
  \label{eq:bootstrap-equivalence}
\end{equation}

\begin{lemma}[Energy inequalities on the bootstrap interval]
\label{lem:bootstrap-evolution}
On the bootstrap interval, the energies satisfy
\begin{align}
  K_m' &\le C(C_1\varepsilon)^2\tau^{-3/2}(K_m+A_m),\qquad m=N-1,N,
    \label{eq:bootstrap-evol-kg}\\
  B_N' &\le C(C_1\varepsilon \tau^{-3/2}+(C_1\varepsilon)^2\tau^{-1})B_N
             +CC_1\varepsilon \tau^{-1/2}K_N,
    \label{eq:bootstrap-evol-wave}\\
  A_{N-1}' &\le CC_1\varepsilon \tau^{-3/2}(A_{N-1}+K_N)
                 +C(C_1\varepsilon)^2\tau^{-2}A_N .
    \label{eq:bootstrap-evol-low}
\end{align}
Here primes may be read as upper right Dini derivatives in \(\tau\).
The undifferentiated potentials satisfy
\begin{equation}
  \|X(\tau)\|_\infty\le C\|X(\tau_0)\|_\infty+C(C_1\varepsilon)^2 .
  \label{eq:bootstrap-potentials}
\end{equation}
\end{lemma}

\begin{proof}
$\bullet$ \emph{Klein--Gordon energy \(K_m\).} In Lemma~\ref{lem:general-kg-energy}, take \(u=\psi\),
\(p=\tau^{-2}\), \(q=1+(4\tau^2)^{-1}\), and \(f=\mathcal R\).
The weights are decreasing and \(E_m=K_m\). Thus
\eqref{eq:general-kg-energy} and Lemma~\ref{lem:kg-remainder}
give \eqref{eq:bootstrap-evol-kg}.

$\bullet$ \emph{Top wave energy \(B_N\).} For the wave system, set \(W_j=\partial_y^jU\) and write
\[
  \partial_\tau W_j=\tau^{-1}M(v)\partial_yW_j+F_j,\qquad
  F_j=\tau^{-1}\big(\partial_y^{j+1}(M(v)U)-M(v)\partial_y^{j+1}U\big).
\]
Every term in \(\tau F_j\) distributes \(j\) derivatives between \(M_y\)
and \(U\). Since \(\|M_y\|_\infty\lesssim\|v_y\|_\infty\),
\(\|M_y\|_{H^m}\lesssim\|v\|_{H^{m+1}}\), and
\(\|U\|_\infty\lesssim C_1\varepsilon\), the tame product estimate gives
\[
  \left(\sum_{j=0}^m\|F_j\|_2^2\right)^{1/2}
  \lesssim \tau^{-1}\big(\|v_y\|_\infty A_m+C_1\varepsilon\|v\|_{H^{m+1}}\big).
\]
Apply Lemma~\ref{lem:general-system-energy} to
\(W=(W_0,\ldots,W_m)\), using the same matrices on each diagonal block.
At \(m=N\), take \(P=\tau^{-1}M(v)\) and \(G=H(v)\) on each block, so \(E_G=B_N\). The symmetry of \(HM\) gives $J=0$, so the term involving \(\|W_y\|_2\) in \eqref{eq:general-system-energy} vanishes. Since \(\|H_\tau\|_\infty\le2\|vv_\tau\|_\infty\) and
\(\|(HM)_y\|_\infty\lesssim\|v_y\|_\infty\), \eqref{eq:general-system-energy}
and the commutator bound give
\[
   B_N'\lesssim
  (\|vv_\tau\|_\infty+\tau^{-1}\|v_y\|_\infty)A_N+C_1\varepsilon\tau^{-1}\|v\|_{H^{N+1}} .
\]
Applying Proposition~\ref{prop:bootstrap-pointwise},
Lemma~\ref{lem:v-sobolev-top}, and \eqref{eq:bootstrap-equivalence} yields \eqref{eq:bootstrap-evol-wave}.

$\bullet$ \emph{Lower-order wave energy \(A_{N-1}\).} At order \(m=N-1\), we apply \eqref{eq:general-system-energy}
with \(G=I\) on each diagonal block. Then \(E_G=A_{N-1}\)
and \(G_\tau=0\). Taking the symmetric and skew-symmetric parts of \(P=\tau^{-1}M(v)\) gives \(\|S_y\|_\infty\lesssim\tau^{-1}\|v_y\|_\infty\) and \(\|J\|_\infty\lesssim\tau^{-1}\|v\|_\infty^2\) by \eqref{eq:potential-skew}.
The remaining derivative satisfies
\(\|W_y\|_2^2=\sum_{j=0}^{N-1}\|\partial_y W_j\|_2^2\le A_N^2\).
Using the commutator bound gives
\[
  A_{N-1}'\lesssim \tau^{-1}\bigl(
     \|v_y\|_\infty A_{N-1}+\|v\|_\infty^2 A_N
     +C_1\varepsilon\|v\|_{H^N}\bigr).
\]
Proposition~\ref{prop:bootstrap-pointwise} and
Lemma~\ref{lem:v-sobolev} with \(m=N\) yield \eqref{eq:bootstrap-evol-low}.

$\bullet$ \emph{Bound for the potentials \(X\).} Write \(X_\tau-\tau^{-1}M(0)X_y=\tau^{-1}h\), where \(h=(M(v)-M(0))U\).
Since \(M(0)\) interchanges the two components, \(Z_\pm=\phi\pm w\)
satisfy
\[
  (\partial_\tau\mp\tau^{-1}\partial_y)Z_\pm=\tau^{-1}h_\pm,\qquad h_\pm=h_1\pm h_2 .
\]
Integrating along the characteristics gives
\[
  Z_\pm(\tau,y)=Z_\pm(\tau_0,y\pm \log(\tau/\tau_0))
    +\int_{\tau_0}^{\tau}h_\pm(\sigma,y\pm\log(\tau/\sigma))\,\frac{d\sigma}{\sigma} .
\]
Translations in \(y\) preserve the \(L^\infty\) norm. Since
\(X=\tfrac12(Z_++Z_-,Z_+-Z_-)^{\mathsf T}\) and
\(\|h\|_\infty\lesssim\|v\|_\infty\|U\|_\infty\) for
\(\|v\|_\infty\le\tfrac12\), taking norms gives
\[
  \|X(\tau)\|_\infty
  \le C\|X(\tau_0)\|_\infty+
    C\int_{\tau_0}^\tau\|v(\sigma)\|_\infty\|U(\sigma)\|_\infty
             \frac{d\sigma}{\sigma}.
\]
By \eqref{eq:pointwise-v} and \eqref{eq:pointwise-U}, the integral is
bounded by \(C(C_1\varepsilon)^2\int_{\tau_0}^{\infty}\sigma^{-3/2}\,d\sigma
\le 2C(C_1\varepsilon)^2\). This proves
\eqref{eq:bootstrap-potentials}.
\end{proof}

\begin{proposition}[Bootstrap improvement]
\label{prop:bootstrap-improvement}
Fix \(\tau_0=2\). Let \(N\ge6\), \(C_0>0\), and \(0<\delta<\tfrac12\).
There exist \(C_1\ge1\) and \(\varepsilon_0>0\) such that, for
\(0<\varepsilon\le\varepsilon_0\), the initial bound
\eqref{eq:bootstrap-data} and the bootstrap assumptions
\eqref{eq:bootstrap-low}--\eqref{eq:bootstrap-kg} imply
\begin{align}
  A_{N-1}(\tau)+K_{N-1}(\tau)+\|X(\tau)\|_\infty
     &\le \tfrac12 C_1\varepsilon, \label{eq:improved-low}\\
  B_N(\tau)&\le \tfrac12 C_1\varepsilon \tau^{1/2+\delta},
     \label{eq:improved-wave}\\
  K_N(\tau)&\le \tfrac12 C_1\varepsilon \tau^\delta
     \label{eq:improved-kg}
\end{align}
throughout \([\tau_0,\tau^*]\), with \(C_1\) and \(\varepsilon_0\) independent of \(\tau^*\).
\end{proposition}

\begin{proof}
By \eqref{eq:bootstrap-equivalence} and the bootstrap assumptions,
\(A_N\lesssim C_1\varepsilon \tau^{1/2+\delta}\), \(K_N\le C_1\varepsilon \tau^\delta\), and
\(A_{N-1}+K_{N-1}\le C_1\varepsilon\).

$\bullet$ \emph{Low-order bounds.}
Substitution in \eqref{eq:bootstrap-evol-low} and
\eqref{eq:bootstrap-evol-kg} gives
\[
  A_{N-1}'\lesssim(C_1^2\varepsilon^2+C_1^3\varepsilon^3)\tau^{-3/2+\delta},
  \qquad
  K_{N-1}'\lesssim(C_1\varepsilon)^3\tau^{-3/2}.
\]
Since \(\delta<\tfrac12\), both right-hand sides are integrable up to
infinity. Integrating from \(\tau_0\) and adding
\eqref{eq:bootstrap-potentials}, we obtain
\begin{equation}
  A_{N-1}(\tau)+K_{N-1}(\tau)+\|X(\tau)\|_\infty
  \le C C_0\varepsilon+C_\delta(C_1^2\varepsilon^2+C_1^3\varepsilon^3).
  \label{eq:close-low-bound}
\end{equation}

$\bullet$ \emph{Top Klein--Gordon energy.}
Taking \(m=N\) in \eqref{eq:bootstrap-evol-kg} yields
\[
  K_N'\lesssim C_1^3\varepsilon^3(\tau^{-3/2+\delta}+\tau^{-1+\delta})
            \lesssim (C_1\varepsilon)^3 \tau^{-1+\delta}.
\]
Since \(\delta>0\), integration gives
\begin{equation}
  K_N(\tau)\le C_0\varepsilon+C_\delta C_1^3\varepsilon^3 \tau^\delta .
  \label{eq:close-kg-bound}
\end{equation}

$\bullet$ \emph{Top wave energy.}
Use \(K_N\le C_1\varepsilon \tau^\delta\) in \eqref{eq:bootstrap-evol-wave} and
choose \(\mu=C(C_1\varepsilon)^2<\delta\), with \(C\) the coefficient of
\((C_1\varepsilon)^2\tau^{-1}B_N\) there. The weighted energy \(Y=\tau^{-\mu}B_N\)
satisfies
\[
  Y'\le CC_1\varepsilon \tau^{-3/2}Y+CC_1^2\varepsilon^2\tau^{-1/2+\delta-\mu}.
\]
Gronwall's inequality gives
\[
  Y(\tau)\le
  \exp\!\left(\int_{\tau_0}^\tau CC_1\varepsilon\sigma^{-3/2}\,d\sigma\right)
  \left(Y(\tau_0)+CC_1^2\varepsilon^2\int_{\tau_0}^\tau
        \sigma^{-1/2+\delta-\mu}\,d\sigma\right).
\]

The exponential factor is uniformly bounded since \(C_1\varepsilon\le1\), \(\tau_0=2\), and
\[
  \int_{\tau_0}^\tau CC_1\varepsilon\sigma^{-3/2}\,d\sigma
  =2CC_1\varepsilon(\tau_0^{-1/2}-\tau^{-1/2})\le 2CC_1\varepsilon \tau_0^{-1/2}.
\]
Using \(Y(\tau_0)=\tau_0^{-\mu}B_N(\tau_0)\le C_0\varepsilon\) and multiplying
by \(\tau^\mu\), we obtain
\begin{align}
  B_N(\tau)
   &\le C C_0\varepsilon \tau^\mu+
    CC_1^2\varepsilon^2\tau^\mu\int_{\tau_0}^\tau\sigma^{-1/2+\delta-\mu}\,d\sigma
    \notag\\
   &\le C C_0\varepsilon \tau^\mu+CC_1^2\varepsilon^2\tau^{1/2+\delta}.
   \label{eq:close-wave-bound}
\end{align}
In the last step we used \(\mu<\delta\) to bound
\[
  \tau^\mu\int_{\tau_0}^\tau\sigma^{-1/2+\delta-\mu}\,d\sigma
  \le \frac{\tau^{1/2+\delta}}{\tfrac12+\delta-\mu}
  \le 2\tau^{1/2+\delta}.
\]

$\bullet$ \emph{Choice of constants.}
By \eqref{eq:close-low-bound}--\eqref{eq:close-wave-bound}, each of
\[
  A_{N-1}+K_{N-1}+\|X\|_\infty,\qquad
  \tau^{-1/2-\delta}B_N,\qquad \tau^{-\delta}K_N
\]
is bounded by \(C_*\varepsilon+C_\delta(C_1^2\varepsilon^2+C_1^3\varepsilon^3)\), where
\(C_*\) and \(C_\delta\) are independent of \(C_1,\varepsilon,\tau^*\).
Choose \(C_1\ge4\max\{1,C_*\}\). With this \(C_1\) fixed, choose
\(\varepsilon_0\) small enough that \(\|v\|_\infty\le\tfrac12\)
follows from \eqref{eq:pointwise-v}, \(C(C_1\varepsilon)^2<\delta\),
and \(C_\delta((C_1\varepsilon)^2+(C_1\varepsilon)^3)\le C_1\varepsilon/4\)
for all \(0<\varepsilon\le\varepsilon_0\). This gives
\eqref{eq:improved-low}--\eqref{eq:improved-kg} simultaneously.
\end{proof}

The initial bound \eqref{eq:bootstrap-data} starts the bootstrap by
local existence and continuity. Proposition~\ref{prop:bootstrap-improvement}
precludes a first time at which any bootstrap bound is saturated.
The estimates therefore hold throughout the lifespan. At any proposed
finite endpoint, they bound \(U\) and \(\psi_\tau\) in \(H^N\),
\(\psi\) in \(H^{N+1}\), and \(X\) in \(L^\infty\). Moreover,
\(X_\tau=\tau^{-1}M(v)U\) and \(|v|\le\tfrac12\), so the first-order system
remains uniformly symmetrizable. The local continuation criterion
then extends the solution past that endpoint. Thus the solution and
the bounds hold for all \(\tau\ge \tau_0\). To conclude the proof of
Theorem~\ref{thm:main-result}, take 
\(\delta=\tfrac18\). The initial hyperbolic data then satisfy
\eqref{eq:bootstrap-data} with \(C_0\) depending only on
\(N\). Thus the smallness threshold depends only on \(N\).
Together with the local Cauchy solution below the initial hyperbola,
this gives the global solution on \([2,\infty)\times\mathbb R\).

\begin{corollary}[Decay estimate for the Klein--Gordon field \(v\)]
\label{cor:physical-time-decay}
Under the assumptions of Theorem~\ref{thm:main-result}, the Klein--Gordon
field satisfies
\[
  \|v(t,\cdot)\|_{L^\infty(\mathbb R_x)}
  \lesssim\varepsilon(1+t)^{-1/2},\qquad t\ge2.
\]
\end{corollary}

\begin{proof}
As in Section~\ref{sec:bootstrap-pointwise}, take \(\tau_0=2\).
By finite speed of propagation, the support of
\(\psi=\tau^{1/2}e^{\phi/2}v\) is contained in the region where
\(\tau e^{-|y|}=t-|x|\ge1\).
For each fixed \(y\), we then have \(\psi(\tau,y)=0\) whenever
\(2\le\tau<e^{|y|}\).

Rewrite \eqref{eq:rescaled-kg} in the form
\[
  \psi_{\tau\tau}+\psi=F,\qquad
  F=\tau^{-2}(\psi_{yy}-\tfrac14\psi)+\mathcal R.
\]

By \eqref{eq:improved-low}, the definition of \(K_{N-1}\), and Sobolev
embedding in \(y\), we have
\[\|\psi\|_\infty+\|\psi_{yy}\|_\infty\lesssim\varepsilon.\]
Lemma~\ref{lem:kg-remainder} with \(m=N-1\) also gives
\[\|\mathcal R\|_\infty\lesssim\varepsilon^3\tau^{-3/2}.\]
Combining these estimates gives
\[
  \|F(\tau,\cdot)\|_\infty
     \lesssim\varepsilon\tau^{-2}+\varepsilon^3\tau^{-3/2}.
\]
If \(e^{|y|}>\tau_0\), smoothness and the vanishing of \(\psi\) for
\(\tau<e^{|y|}\) yield
\[\psi(e^{|y|},y)=\psi_\tau(e^{|y|},y)=0.\] 
For \(\tau\ge e^{|y|}\), Duhamel's
formula at fixed \(y\) gives
\[
  |\psi(\tau,y)|
  =\left|\int_{e^{|y|}}^\tau \sin(\tau-\sigma)F(\sigma,y)\,d\sigma\right|
  \lesssim \varepsilon e^{-|y|}+\varepsilon^3e^{-|y|/2}
  \lesssim \varepsilon e^{-|y|/2}.
\]
If \(e^{|y|}\le\tau_0\), then \(|y|\le\log 2\), and the uniform
bound on \(\psi\) gives the same estimate. Using the uniform bound on
\(\phi\) from \eqref{eq:improved-low} and \(\tau e^{|y|}=t+|x|\), we obtain
\[
|v(t,x)|=\tau^{-1/2}e^{-\phi/2}|\psi(\tau,y)|
  \lesssim \varepsilon (t+|x|)^{-1/2}
  \lesssim \varepsilon (1+t)^{-1/2}.
\]
Finally, on the portion of the support with \(\tau<2\), we have
\(2t-1\le t^2-x^2<4\), so \(2\le t<\tfrac52\). On this bounded time
interval, the local estimate in Section~\ref{sec:bootstrap-pointwise}
gives 
\[\|v(t,\cdot)\|_\infty\lesssim\varepsilon\lesssim\varepsilon(1+t)^{-1/2}.\]
\end{proof}


\section*{Acknowledgements} 
The author would like to express her sincere gratitude to her postdoctoral advisor, Siyuan Ma, for continuous support, encouragement, and valuable guidance during her postdoctoral research. The author is also deeply grateful to her Ph.D. advisor, Shijie Dong, for constant encouragement, insightful advice, and generous help over the years.

\section*{Data Availability}
The data are available upon request.


\end{document}